\documentclass[a4paper,12pt]{amsart}

\usepackage{amsmath}
\usepackage[T1]{fontenc}
\usepackage[francais]{babel}
\def\ni{\noindent}
\def\vs{\vskip .6cm}
\def\ss{\smallskip}
\def\n{\nabla}
\def\l{\lambda}

\def\e{\epsilon}
\def\beq{\begin{eqnarray*}}
\def\eeq{\end{eqnarray*}}

\def\x{\times}

\def\r{\end{proof}}

\def \RM{\mathbb{R}}

\def \ZM{\mathbb{Z}}
\def \CM{\mathbb{C}}

\def \HM{\mathbb{H}}

\def \m{\mathfrak m}
\def \h{\mathfrak h}
\def \g{\mathfrak g}
\def \su {\mathfrak{su}}
\def \u {\mathfrak{u}}
\def \n {\mathfrak{n}}
\def \l {\mathfrak{l}}
\def \k {\mathfrak{k}}

\def \la {\langle}
\def \ra {\rangle}

\newtheorem{defi}{D{\'e}finition}[section]

\newtheorem{prop}[defi]{Proposition}

\newtheorem{theo}[defi]{Theor{\`e}me}

\newtheorem{lemm}[defi]{Lemme}

\newtheorem{conj}[defi]{Conjecture}
 
\newenvironment{demo}{{\it D{\'e}monstration.}}{$\Box$ \ \\}

\begin{document}

\begin{center}
{\large\bf Classification des vari{\'e}t{\'e}s approximativement
  k{\"a}hleriennes homog{\`e}nes} \vs

{\bf Jean-Baptiste Butruille}\ss

{\footnotesize
Centre de Math{\'e}matiques, {\'E}cole Polytechnique, UMR 7640 du CNRS, 91128
Palaiseau  
\vskip-0.5ex
e-mail : jbbutruille@math.polytechnique.fr}\vs 
\end{center}

{\footnotesize\sl 

\ni{\bf R{\'e}sum{\'e}} -- On d{\'e}montre la conjecture de Gray et Wolf que
les seules vari{\'e}t{\'e}s strictement approximativement k{\"a}hleriennes
homog{\`e}nes sont les espaces 3-symm{\'e}triques. Pour cela on les
classifie en dimension 6 puis la d{\'e}monstration pour les dimensions
sup{\'e}rieures provient d'un th{\'e}or{\`e}me de Nagy, s'appuyant sur des
r{\'e}sulats pr{\'e}c{\'e}dents de Cleyton et Swann. Les
espaces homog{\`e}nes de dimension 6, $S^3 \x S^3$, $S^6$, $\CM P(3)$ et
l'espace des drapeaux $F(1,2)$ portent une unique structure approximativement
k{\"a}hlerienne invariante {\`a} homot{\'e}tie pr{\`e}s. Pour le premier, cela
r{\'e}sulte de la r{\'e}solution d'une {\'e}quation diff{\'e}rentielle
donn{\'e}e par Reyes Carri{\'o}n. Pour les deux derniers, il s'agit de leur
structure presque hermitienne d'espace de twisteurs sur une
vari{\'e}t{\'e} de dimension 4. Enfin, les structures
approximativement k{\"a}hleriennes sur la sp{\`e}re de dimension 6
correspondent {\`a} des 3-formes constantes sur $\RM^7$.}\vs \vs

\section{Introduction}

Une vari{\'e}t{\'e} approximativement k{\"a}hlerienne (ou Nearly K{\"a}hler, en
anglais : dor{\'e}navant on note NK) est une vari{\'e}t{\'e} presque hermitienne 
telle que la deriv{\'e}e covariante pour la connexion de Levi-Civit{\'a} de la forme
de K{\"a}hler $\omega$ est antisym{\'e}trique.

Les vari{\'e}t{\'e}s strictement approximativement k{\"a}hleriennes (SNK) -- c'est
{\`a} dire qui ne sont pas simplement k{\"a}hleriennes -- en dimension 6 sont 
particuli{\`e}res {\`a} plusieurs titres. Notamment elles admettent un spineur
de Killing (voir \cite{gru}) et sont d'Einstein (\cite{gr2}). Depuis
les travaux de Nagy dans \cite{na} on sait que de 
leur classification d{\'e}pend beaucoup celle des vari{\'e}t{\'e}s NK en toute dimension. 

Une voie privil{\'e}gi{\'e}e de construction des vari{\'e}t{\'e}s NK est celle
des espaces 3-sym{\'e}triques. A. Gray a montr{\'e} en 1972 que tout espace
homog{\`e}ne 3-sym{\'e}trique naturellement r{\'e}ductif est muni
canoniquement d'une structure NK (voir \cite{gr3}). Auparavant A. Gray et
J.A. Wolf avaient conjectur{\'e} en 1968 dans \cite{wo} que tout espace homog{\`e}ne
SNK est 3-sym{\'e}trique. Retournant a l'int{\'e}r{\^e}t manifest{\'e}
dans les ann{\'e}es 70 pour ces vari{\'e}t{\'e}s, P.A. Nagy (\cite{na2}) a
d{\'e}compos{\'e} les vari{\'e}t{\'e}s SNK de dimension quelconque en produits
riemanniens d'espaces de twisteurs au-dessus de vari{\'e}t{\'e}s
K{\"a}hler-quaternioniques, d'espaces homog{\`e}nes et de
vari{\'e}t{\'e}s SNK de dimension 6. Si la vari{\'e}t{\'e} est de plus homog{\`e}ne,
il ne reste que des espaces homog{\`e}nes, de divers types mais tous
3-sym{\'e}triques et des vari{\'e}t{\'e}s SNK homog{\`e}nes de dimension 6 de
telle sorte qu'il suffit pour que la conjecture soit vraie qu'elle le
soit en dimension 6.

Dans cet article on prouve la conjecture en classifiant
totalement les vari{\'e}t{\'e}s SNK homog{\`e}nes simplement connexes de
dimension 6.

\begin{theo}
Tous les espaces homog{\`e}nes strictement NK de dimension 6 sont des
espaces 3-sym{\'e}triques munis de leur structure presque complexe canonique.
\label{theo}
\end{theo}

\begin{theo}
Les seuls espaces homog{\`e}nes strictement NK simplement connexes
de dimension 6 sont isomorphes {\`a} $G/H$ o{\`u} $G$ et $H$ sont les
groupes de Lie donn{\'e}s dans la liste :
\begin{itemize}
\item $G=S^3 \x S^3$ et $H=\{1\}$
\item $G=G_2$ et $H$  est $SU(3)$ (dans ce cas $G/H$ est la sph{\`e}re
  de dimension 6) ou un de ses sous-groupes finis
\item $G=Sp(2)$ et $H=S^1 \x SU(2)$ (alors $G/H$ est l'espace
  projectif complexe $\CM P(3)$) ou un de ses sous-groupes finis
\item $G=SU(3)$, $H=S^1 \x S^1$ et $G/H$ est l'espace de drapeaux
  $F(1,2)$
\end{itemize}
De plus sur chacun de ces espaces homog{\`e}nes il y a une seule
  structure presque complexe NK invariante {\`a} isomorphisme pr{\`e}s.
\label{theo2}
\end{theo}

Pour achever la classification il faudra examiner les quotients finis
des espaces homog{\`e}nes dont la liste est dans le th{\'e}or{\`e}me
\ref{theo2}. En effet les vari{\'e}t{\'e}s SNK sont compactes, de groupe
fondamental fini d'apr{\`e}s Nagy (\cite{na}). En dimension 6 cela
r{\'e}sultait d{\'e}j{\`a} de la d{\'e}monstration par Gray dans \cite{gr2}
qu'elles sont d'Einstein, {\`a} courbure scalaire strictement positive.

Dans la section 2, des pr{\'e}liminaires alg{\'e}briques {\`a} la
classification permettent d'{\'e}tablir une premi{\`e}re liste des groupes
$G$ et $H$ tels que leur quotient $G/H$ est suceptible d'admettre une
structure SNK invariante. Dans la section 3, on traite le cas jug{\'e} le plus
difficile du produit de sph{\`e}res $S^3 \x S^3$. En r{\'e}solvant
l'{\'e}quation diff{\'e}rentielle de Reyes Carri{\'o}n \cite{re}, qui caract{\'e}rise les
vari{\'e}t{\'e}s SNK en dimension 6, sur l'espace des 2-formes
invariantes, on montre que $S^3 \x S^3$ admet une seule structure
homog{\`e}ne SNK correspondant {\`a} la construction de Ledger et Obata dans
\cite{le} d'un espace 3-sym{\'e}trique. Dans la section 4 on traite plusieurs cas
qui se ram{\`e}nent au pr{\'e}c{\'e}dent. Dans la section 5 et la section 6,
le moyen de d{\'e}terminer les structures presques hermitiennes
invariantes de $F(1,2)$ et $\CM P(3)$ est de d{\'e}composer la
repr{\'e}sentation lin{\'e}aire isotropique en repr{\'e}sentations
irr{\'e}ductibles. Cette recherche appara{\^\i}t
li{\'e}e {\`a} la th{\'e}orie des espaces de twisteurs des vari{\'e}t{\'e}s de
dimension 4 (ici $\CM P(2)$ et $S^4$ respectivement.) Il faut alors
calculer quelles sont NK. Enfin dans la section 7, la donn{\'e}e d'une
structure NK de la sph{\`e}re $S^6$ est rappel{\'e}e {\^e}tre {\'e}quivalente
{\`a} la donn{\'e}e d'une 3-forme g{\'e}n{\'e}rique constante sur l'espace euclidien
$\RM^7$ et cela termine notre {\'e}tude.  

\section{Pr{\'e}liminaires}

Une vari{\'e}t{\'e} riemannienne $(M^n,g,J)$ est dite presque hermitienne si la
m{\'e}trique $g$ et la structure presque complexe $J$ v{\'e}rifient
$g(JX,JY)=g(X,Y)$ quels que soient les champs de vecteurs $X$ et
$Y$, donnant une reduction du fibr{\'e} principal $SO(M)$ {\`a} $U(n)$,
not{\'e}e $U(M)$. Soit $\nabla$ la connexion de Levi-Civit{\'a} de $g$,
$\omega=g(J.,.)$ la forme de K{\"a}hler. Soit $\Lambda^{(2,0)+(0,2)} M$ le
sous-fibr{\'e} des 2-formes $\nu$ de $M$ v{\'e}rifiant
$\nu(JX,JY)=-\nu(X,Y)$. Le tenseur $\nabla \omega$
est une section de $T^*M \otimes \Lambda^{(2,0)+(0,2)}M$. Il s'agit d'un fibr{\'e}
associ{\'e} de $U(M)$ et on peut voir alternativement $\nabla \omega$ comme
une fonction {\'e}quivariante de $U(M)$ dans $V^* \otimes
\Lambda^{(2,0)+(0,2)} V^*$ o{\`u} $V$ est un espace vectoriel hermitien
avec $dim_{\RM} V=n$. Or il se d{\'e}compose, comme espace de
repr{\'e}sentation de $U(n)$, en quatre sous-espaces irr{\'e}ductibles non
isomorphes, not{\'e}s $W_1,W_2,W_3,W_4$. Selon que $\nabla \omega$ prend
ses valeurs dans un des 16 sous-espaces invariants $\bigoplus_{i \in
  I} W_i$, $I \subset \{1,2,3,4\} $ on d{\'e}finit, apr{\`e}s Gray et
Hervella \cite{gr3}, 16 classes naturelles de vari{\'e}t{\'e}s presque
hermitiennes. En reprenant leurs notations, les vari{\'e}t{\'e}s NK sont
d{\'e}finies par $\nabla \omega \in W_1 = \{\nu \in V^* \otimes
\Lambda^{(2,0)+(0,2)} V^* | \nu(X,Y,Z)=-\nu(Y,X,Z) \}$. Autrement dit le
tenseur $\nabla J$ de type (2,1) est totalement antisym{\'e}trique. Maintenant,
par l'interm{\'e}diaire de fonctions $U(n)$-{\'e}quivariantes liant
$\nabla \omega$ {\`a} d'autres tenseurs g{\'e}om{\'e}triquement significatifs de la
vari{\'e}t{\'e}, on se rend compte que, par exemple, la composante
sur $W_1 \oplus W_2$ est donn{\'e}e par le tenseur de Nijenhuis ou la
composante sur $W_3 \oplus W_4$ par la partie de type (2,1)+(1,2) de
$d \omega$. On d{\'e}montre de cette fa\c con que la structure presque
complexe d'une vari{\'e}t{\'e} NK n'est jamais int{\'e}grable {\`a} moins
qu'elle soit k{\"a}hlerienne (i.e. avec $\nabla \omega =0$) et que
la diff{\'e}rentielle de la forme de K{\"a}hler $d \omega$ est de type
(3,0)+(0,3).

\vs
Soit $M=G/H$ un espace homog{\`e}ne r{\'e}ductif c'est {\`a} dire qu'il
existe une d{\'e}composition $\g=\h \oplus \m$, $Ad(H)$-invariante. La
projection naturelle $\pi : G \to G/H$ est une fibration principale de
groupe $H$. La restriction de la diff{\'e}rentielle de la projection en
$e$ est un isomorphisme $\pi_* : \m \to T_oM$ o{\`u} $o=\pi(e)$. Cette
identification est de plus $H$-{\'e}quivariante : $\pi_*(Ad_h
u)=h_*\pi_*(u), \ \forall \ u \in \m$ (en notant encore $g$ le
diff{\'e}omorphisme de $M$ induit par la multiplication dans le groupe {\`a}
gauche par $g \in G$.) Ainsi la repr{\'e}sentation lin{\'e}aire
isotropique est vue comme la repr{\'e}sentation $Ad(H)$ sur
$\m$. Autrement dit le fibr{\'e} associ{\'e} $G \times_{Ad} \m/H$ est
canoniquement isomorphe au fibr{\'e} tangent par l'application
\begin{eqnarray*}
G \times_{Ad} \m/H    & \stackrel{\sim}{\to} & TM \\
\lbrack g,u \rbrack & \mapsto              & g_*\pi_*(u)
\end{eqnarray*}
Une section d'un espace de tenseurs de $M$ est identifi{\'e}e {\`a}
une fonction $H$-{\'e}quivariante de $G$ sur l'espace de tenseurs
correspondant de $\m$, c'est {\`a} dire d'un espace vectoriel fixe. Si
cette section est de plus invariante pour l'action de $G$ sur $M$
induite par la multiplication {\`a} gauche, la fonction est constante
{\'e}gale {\`a} un tenseur $Ad(H)$-invariant.

\vs
A l'intersection de ces deux notions, un espace homog{\`e}ne presque
hermitien $(M,g,J)$ est un espace homog{\`e}ne dont la m{\'e}trique et la
structure presque complexe sont invariantes pour l'action de $G$ et en
font une vari{\'e}t{\'e} presque hermitienne.

Elles sont donc identifi{\'e}es {\`a} un produit scalaire $(.|.)$ et un
endomorphisme de carr{\'e} $-1$ de $\m$, $Ad(H)$-invariants. La
d{\'e}riv{\'e}e covariante $\tilde{\nabla}$ d'une connexion n'est pas un
tenseur, en revanche $A$ d{\'e}fini par
$A_X=\mathcal{L}_X-\tilde{\nabla}_X$, si. Il y a donc une
correspondance bijective entre les connexions $G$-invariantes et les
applications $\Lambda : \m \to \mathfrak{so}(\m)$, c'est {\`a} dire
v{\'e}rifiant
\begin{equation}
(\Lambda(X)Y | Z) + (\Lambda(X)Z | Y) = 0
\quad \text{pour tous } X,Y \in \m
\label{L'}
\end{equation}
donn{\'e}e par $\Lambda(X)=(-A_X)_o$, en identifiant $\m$ et $T_o M$. Il
s'agit du th{\'e}or{\`e}me de Wang sp{\'e}cialis{\'e} aux espaces homog{\`e}nes
r{\'e}ductifs (th{\'e}or{\`e}me 2.1, p191 de \cite{ko}). Si
$\tilde{\nabla}=\nabla$, la d{\'e}riv{\'e}e covariante de
la connexion de Levi-Civit{\'a}, qui est sans torsion, on a de plus
\begin{equation}
\Lambda(X)Y - \Lambda(Y)X = [X,Y]_\m \quad \text{pour tous } X,Y \in \m
\label{L}
\end{equation}
La notation $\m$ en indice signifie la projection sur ce sous-espace
du vecteur indic{\'e}.
Les {\'e}quations (\ref{L'}) et (\ref{L}) d{\'e}finissent une unique
application $\Lambda$ donn{\'e}e par la formule du th{\'e}or{\`e}me 3.3, p201 \cite{ko}
\begin{eqnarray}
\Lambda(X)Y = \frac{1}{2}[X,Y]_{\m} + U(X,Y) \label{LC} \\
2(U(X,Y)|Z)=([Z,Y]_{\m}|X)+([Z,X]_{\m}|Y)
\label{U}
\end{eqnarray}

\vs
Si $G$ est muni d'un endomorphisme $\sigma$ d'ordre 3 dont $H$ est l'ensemble
des points fixes et que le produit scalaire $Ad(H)$-invariant sur $\m$
repr{\'e}sentant la m{\'e}trique $(.|.)$ est de plus $\sigma_*$-invariant,
en passant au quotient on obtient une isom{\'e}trie $\theta$ v{\'e}rifiant
$\theta^3=Id$ dont $o$ est un point fixe isol{\'e}. En effet
$\theta(\pi(g))=\pi(g)$ implique $\sigma(g)=hg$ o{\`u} $h \in H$ est une
racine cubique de l'unit{\'e}. Par cons{\'e}quent si on est assez proche
de $e$, $h=e$ et $g$ doit appartenir {\`a} $H$. En conjuguant $\theta$
par l'action du groupe, transitive, chaque point de $M$ est rapport{\'e}
{\`a} une telle application. Cela
signifie que $M$ est un espace 3-sym{\'e}trique au sens de la
\begin{defi}
Un espace 3-sym{\'e}trique est une vari{\'e}t{\'e} $M$ munie d'une famille
d'isom{\'e}tries globales $\theta_m$, $m \in M$ telles que $m$ est un
point fixe isol{\'e} de $\theta_m$ et $\forall \ m \in M$, $\theta_m^3=Id$.
\end{defi}

Comme on {\'e}crit dans le plan complexe une racine cubique de l'unit{\'e}, on
associe {\`a} tout espace 3-sym{\'e}trique une structure presque complexe
dite canonique en posant
\[
\forall \ m \in M \quad (\theta_m)_*=-\frac{1}{2}Id_{T_m M} + 
\frac{\sqrt{3}}{2}J_m
\]
En notant encore $J$ l'endomorphisme de carr{\'e} -1 de $\m$ associ{\'e}, Gray
a calcul{\'e} dans \cite{gr} que
\begin{equation}
\forall \ X,Y \in \m, \quad [X,JY]_{\m}=-J[X,Y]_{\m}
\label{T}
\end{equation}

Il en a d{\'e}duit
\begin{theo}
Soit $M$ un espace homog{\`e}ne 3-sym{\'e}trique. Les propositions suivantes
sont {\'e}quivalentes. \\
(\romannumeral 1) la structure presque complexe canonique de $M$ est
NK \\
(\romannumeral 2) $M$ est naturellement r{\'e}ductif \\
\end{theo}

\begin{demo}
Par naturellement r{\'e}ductif on entend que la m{\'e}trique satisfait
\begin{equation}
([X,Y]_{\m}|Z)=(X|[Y,Z]_{\m}) \quad \text{pour tous } X,Y,Z \in \m
\end{equation}
C'est le cas en particulier si $(.|.)$ est la restriction {\`a} $\m$
d'une m{\'e}trique biinvariante de $G$. A l'aide de $(\ref{LC})$ et
$(\ref{T})$ on peut calculer que pour tous $X,Y \in \m$ $(\nabla_X J)JY
= 2U(X,Y) + T(X,Y)$. Par cons{\'e}quent $(\romannumeral 1)$ et
$(\romannumeral 2)$ sont {\'e}quivalents {\`a} $U=0$.
\end{demo}

R{\'e}ciproquement on reformule g{\'e}om{\'e}triquement la conjecture de
Gray et Wolf, {\'e}nonc{\'e}e dans \cite{wo} ({\`a} la fin, p113 ou dans
l'introduction, p79) dans le langage de la th{\'e}orie des
groupes de Lie :
\begin{conj}
Tout espace homog{\`e}ne SNK est un espace 3-sym{\'e}trique naturellement
r{\'e}ductif muni de sa structure presque complexe canonique.
\label{conj}
\end{conj}

\vs

D{\'e}sormais $M=G/H$ d{\'e}signe un espace homog{\`e}ne r{\'e}ductif simplement connexe de
dimension 6. On cherche {\`a} quelles conditions il admet une structure
presque hermitienne NK.
Or
\begin{prop}
Soit $M=G/H$ un espace homog{\`e}ne riemannien de dimension 6, non
isom{\'e}trique {\`a} la sph{\`e}re standard $S^6$. S'il admet une structure
presque complexe NK, elle est unique et invariante pour l'action du groupe $G$.
\label{spineurs}
\end{prop} 
\begin{proof}
L'expos{\'e} ci-dessous suit le livre \cite{bau},
surtout les section 5.2 et 5.3.

En dimension 6, les vari{\'e}t{\'e}s NK ont sont caract{\'e}ris{\'e}es par
l'existence d'un spineur de Killing : si on note $\Sigma M$ le fibr{\'e} des
spineurs complexes, il existe $\psi \in
\Gamma(\Sigma M)$, $\beta \in \RM \backslash \{0\}$ tel que
\[ \nabla_X \psi = \beta X.\psi
\]
en notant encore $\nabla$ la connexion induite sur $\Gamma(\Sigma)$
par la connexion de Levi-Civit{\'a}. Quant au point, il d{\'e}signe la
multiplication de Clifford. Comme on est en dimension paire on peut
scinder $\psi$ en une part n{\'e}gative $\psi^- \in \Sigma^- M$ et une
part positive $\psi^+ \in \Sigma^+ M$ suivant la d{\'e}composition en
sous-espaces irr{\'e}ductibles de la repr{\'e}sentation de Spin(6). Alors
\begin{eqnarray*}
\nabla_X \psi^+ & = & \beta X.\psi^- \\
\nabla_X \psi^- & = & \beta X.\psi^+
\end{eqnarray*}
si bien que le conjugu{\'e} de $\psi$, $\overline \psi=\psi^+ - \psi^-$
est encore un spineur de Killing :
\[ \nabla_X \overline \psi = -\beta X.\overline \psi
\]
De plus $\beta$ et $-\beta$ sont les seules valeurs possibles car il
faut que la courbure scalaire vaille $s=4\beta^2 n(n-1)$, o{\`u} $n=6$ est
la dimension de $M$. On d{\'e}finit alors la structure presque complexe
$J$ associ{\'e}e {\`a} $\psi$ par
\begin{equation}
JX.\psi^+ \stackrel{def}{=} iX.\psi^+
\label{Jpsi}
\end{equation}
apr{\`e}s avoir v{\'e}rifi{\'e} que l'ensemble $\{ X.\psi^+_x | X \in T_xM
\}$ est un sous-espace complexe de $\Sigma^+_x M$ en tout point $x \in
M$. A pr{\'e}sent si $M$ n'est pas la sph{\`e}re $S^6$, l'ensemble des
spineurs de Killing pour la valeur $\beta$ est de dimension 1
(proposition 1, p126) et la structure presque complexe associ{\'e}e est
NK. Inversement si on se donne une structure
presque complexe NK, $J$ sur $M$, il existe un spineur de Killing tel
qu'elle en soit la structure presque complexe associ{\'e}e (voir \cite{gru}).
On en d{\'e}duit qu'il y a une seule structure presque complexe NK sur
$M$. De plus on peut d{\'e}finir une action du groupe d'isom{\'e}trie sur les
spineurs. Pour un espace homog{\`e}ne, les spineurs de Killing sont
invariants et par cons{\'e}quent $J$ aussi, par (\ref{Jpsi}).
\end{proof}
Ainsi, tout espace homog{\`e}ne riemannien de
dimension 6 hormis la sph{\`e}re $S^6$, muni d'une structure presque
complexe NK est un espace homog{\`e}ne presque hermitien.

Dans l'examen de la conjecture en dimension 6 on distingue
deux types d'espaces homog{\`e}nes. Les cas o{\`u} $G$ et $H$ sont des groupes
produits des sph{\`e}res $S^1$ et $S^3$. On veut d{\'e}montrer qu'ils se
ram{\`e}nent tous au cas de $S^3 \x S^3$ c'est {\`a} dire, comme on le verra
{\`a} la section suivante, que la vari{\'e}t{\'e} admet
au plus une structure SNK homog{\`e}ne qui est de plus 3-sym{\'e}trique comme
on souhaite. Et les cas exceptionnels de l'espace des drapeaux
$F(1,2)$, de l'espace projectif $\CM P(3)$ et de la sph{\`e}re $S^6$ qu'on
sait admettre eux aussi une telle structure et on voudrait d{\'e}montrer
en examinant la repr{\'e}sentation lin{\'e}aire isotropique de $H$ qu'il n'y
en n'a pas d'autre.

Une telle d{\'e}marche est syst{\'e}matique : 
\begin{prop}
Soit $G/H$ un espace homog{\`e}ne SNK de dimension 6. Le groupe $H$ est
contenu dans $SU(3)$.
\label{su3}
\end{prop}

\begin{proof}  
D'abord $M=G/H$ est un espace homog{\`e}ne presque hermitien, $\m$ est
muni d'un endomorphisme de carr{\'e} $-1$ et on le voit comme un espace
vectoriel complexe de dimension 3. Alors $Ad(H) \subset U(\m)$, le
groupe des transformations unitaires de $\m$, exprimant que la m{\'e}trique et la
structure presque complexe sont invariantes, ce qui s'{\'e}crit encore
$g^\star \omega = \omega$, o{\`u} $\omega$ est la forme de K{\"a}hler. Puis 
\[ g^\star d \omega = d (g^\star \omega) = d \omega \]
Or pour une vari{\'e}t{\'e} SNK de dimension 6, $d\omega$ est non nulle, de
type $(3,0) + (0,3)$, $Ad(H)$ pr{\'e}serve aussi dans ce cas une 3-forme
complexe sur $\m$ et doit finalement {\^e}tre contenu dans $SU(\m)$. 
\end{proof}
Celui-ci {\'e}tant compact, de
dimension 8, les seuls groupes qui conviennent, hormis lui-m{\^e}me,
sont $U(1)=S^1$, $SU(2)=Sp_1=S^3$ et leurs produits directs
et leurs quotients finis.

Pour un espace homog{\`e}ne on a la suite d'homotopie
\begin{equation}
\pi_2(G/H) \to \pi_1(H) \to \pi_1(G) \to \pi_1(G/H) \to H/H^0 \to 0
\label{homotopie}
\end{equation}
Cela implique que le rang du groupe fondamental de $G$ doit
{\^e}tre inf{\'e}rieur ou {\'e}gal {\`a} celui du groupe fondamental de $H$. Et si on
suppose $G$ connexe, $H$ l'est aussi.

A partir de ces consid{\'e}rations on
peut {\'e}tablir une liste des groupes $H$ possibles et la liste, en
regard, des couples $(G,H)$, compatible avec ces faits. Il nous
reviendra d'examiner cas par cas {\`a} partir de cette liste ce
que l'existence de cette application surjective 
\begin{equation}
\phi : \pi_1(H) \to \pi_1(G)
\label{istar}
\end{equation}
impose plus pr{\'e}cisement au plongement de $H$ dans $G$. Pour {\'e}viter
de citer tous les quotients finis d'un groupe on {\'e}crit seulement la
liste des alg{\`e}bres de Lie :

\begin{lemm}
Soient $G/H$ un espace homog{\`e}ne SNK simplement connexe de dimension
6 et soient $\g$ et $\h$ les alg{\`e}bres de Lie de $G$ et $H$,
respectivement. Elles apparaissent {\`a} la m{\^e}me ligne du
tableau ci-dessous :
\[
\begin{array}{|c|l|l|}
\hline
dim \ \h & \h & \g \\
\hline \hline
0 & \{0\} & \su(2) \oplus \su(2) \\
\hline
1 & i\RM & i\RM \oplus \su(2) \oplus \su(2) \\
\hline
2 & i\RM \oplus i\RM & i\RM \oplus i\RM \oplus \su(2) \oplus \su(2) \\
  & i\RM \oplus i\RM & \su(3) \\ 
\hline
3 & \su(2) & \su(2) \oplus \su(2) \oplus \su(2) \\
\hline
4 & i\RM \oplus \su(2) & i\RM \oplus \su(2) \oplus \su(2) \oplus
\su(2) \\
  & i\RM \oplus \su(2) & \mathfrak{sp}(2) \\
\hline
8 & \su(3) & \g_2 \\
\hline
\end{array}
\]
\label{liste}
\end{lemm}
\begin{demo}
Cette liste est {\'e}tablie, en commen\c cant par $\h$ jusqu'{\`a}
$\su(3)$ qui doit la contenir {\`a} cause de la proposition \ref{su3}, {\`a}
partir de la liste par ordre de dimension croissante des groupes
simples, compacts, connexes : $S^3 \simeq SU(2) \simeq
Sp(1), \ SU(3), \ Spin(5) \simeq Sp(2), \ G_2, \text{ etc.}$ Le
dernier groupe cit{\'e} est de dimension 14 car $SU(3)$ est de
dimension 8 et $G/H$ de dimension 6.

Outre les raisons de dimension qui peuvent l'emp{\^e}cher, on se demande
quelles $\h$ sont vraiment des sous-alg{\`e}bres de Lie de
$\su(3)$. C'est bien s{\^u}r le cas de $i\RM$, $\su(2)$ et $i\RM \oplus
i\RM$. C'est encore le cas de $i\RM \oplus \su(2)$ via le plongement

\begin{eqnarray*}
i\RM \oplus \su(2) & \to     & \su(3) \\
ix + A             & \mapsto & 
\left( \begin{array}{cc}
ix & 0 \\
0  & A - ixI_2
\end{array} \right)
\end{eqnarray*}
o{\`u} $A$ est une matrice $2 \x 2$ antihermitienne {\`a} trace nulle.
En revanche ce n'est plus le cas de $\su(2) \oplus \su(2)$. En effet
{\`a} partir d'une application  $\varphi : \su(2) \oplus \su(2) \to
\su(3)$, non identiquement nulle,
en restreignant {\`a} chaque facteur on obtient deux repr{\'e}sentations
complexes $\rho_1, \rho_2$ de $\su(2)$ de dimension 3 qui
commutent. Alors de deux choses l'une : ou bien une repr{\'e}sentation
est irr{\'e}ductible et la seconde est triviale, par le
lemme de Schur. Ou bien il existe pour chacune un sous-espace de
de dimension 1 et un sous-espace de dimension 2 invariants,
orthogonaux. Cette d{\'e}composition est la m{\^e}me pour $\rho_1$ et
$\rho_2$ car elles commutent et les deux repr{\'e}sentations sont nulles
sur le premier espace car $\su(2)$ n'admet pas de repr{\'e}sentation
autre que triviale avant la dimension 2. Par cons{\'e}quent $\varphi$
est le prolongement par z{\'e}ro d'une application $\su(2) \oplus \su(2) \to
\su(2)$. Dans ce cas, comme dans le pr{\'e}c{\'e}dent, elle ne saurait {\^e}tre
un plongement.
Enfin $i\RM \oplus i\RM \oplus i\RM$ et toutes les alg{\`e}bres qui le
contiennent (en premier $i\RM \oplus i\RM \oplus \su(2)$) ne peuvent
pas {\^e}tre plong{\'e}es dans $\su(3)$ car les sous-alg{\`e}bres de Cartan de
celle-ci sont de dimension 2.
\end{demo}

Si $G$ n'est pas simplement connexe, soit $\pi : \tilde{G} \to G$ son
rev{\^e}tement universel et $H'$ le sous-groupe de $\tilde{G}$,
$\pi^{-1}(H)$. Alors

\begin{lemm}
Les deux espaces homog{\`e}nes $M=G/H$ et $M'=\tilde{G}/H'$ sont
isomorphes.
\label{Gtilde}
\end{lemm}

\begin{proof}
L'isomorphisme est donn{\'e} par
\begin{eqnarray*}
M'                & \to     & M \\
\lbrack g \rbrack & \mapsto & [\pi(g)]
\end{eqnarray*}
qui est bien d{\'e}finie et injective car $g$ et $g'$ d{\'e}finissent la m{\^e}me
classe, $g'g^{-1} \in H'$, si et seulement si $\pi(g')\pi(g)^{-1} \in
H$, c'est {\`a} dire $\pi(g)$, $\pi(g')$ d{\'e}finissent la m{\^e}me classe
de $M$. Elle est aussi surjective car $\pi$ l'est.
\end{proof} 

Si le groupe fondamental de $G$ est fini, $\tilde{G}$ et $H'$ sont
encore compactes, l'alg{\`e}bre de Lie de $H'$, $\h'=\h$ est encore
incluse dans $\su(3)$ et bien s{\^u}r l'alg{\`e}bre de
Lie de $\tilde{G}$ est $\g$. Dans la suite on supposera donc que $G$
est simplement connexe (sections 3, 5, 6, 7) ou $G$ et $H$ sont des
produits finis des groupes $S^1$ et $S^3$ ou des quotients finis
de ces produits (section 4). Dans le premier cas, comme les groupes sont
connexes, la repr{\'e}sentation lin{\'e}aire isotropique de $M$, $Ad(H)$, est
donn{\'e}e par la repr{\'e}sentation de l'alg{\`e}bre de Lie, $ad(\h)$, sur
le suppl{\'e}mentaire invariant choisi. On peut donc finalement se
contenter de chercher les structures presque hermitiennes invariantes
des espaces $G/H$ o{\`u} $G$ et $H$ sont les groupes simples, compacts,
connexes, dont les alg{\`e}bres de Lie apparaissent {\`a} une m{\^e}me ligne
du tableau \ref{liste}.

\section{Le groupe de Lie $S^3 \x S^3$}

On se propose ici de chercher toutes les structures NK sur 
le groupe de Lie $S^3\times S^3$, invariantes {\`a} gauche. On en connait 
d'avance une, correspondant {\`a} la construction de Ledger et Obata (voir
\cite{le}).
 
Si $G$ est un groupe de Lie compact, $\mathfrak{g}$ son alg{\`e}bre de
Lie, $g$ une m{\'e}trique biinvariante sur $G$, on appelle 
$\Delta$ le sous-groupe diagonal, isomorphe a G, de de $G
\times G \times G$. L'espace homog{\`e}ne $M = G \times G \times G/\Delta$ 
est isomorphe {\`a} $G \times G$ et on choisit l'identification
concr{\`e}te
\begin{eqnarray*}
G \x G \x G/\Delta    & \to     & G \x G \\ 
\lbrack x,y,1 \rbrack & \mapsto & (x,y)
\end{eqnarray*}
o{\`u} $[x,y,1]$ d{\'e}signe la classe de $(x,y,1)$. Toute classe contient
un triplet de cette forme donc l'application est bien
d{\'e}finie. On choisit alors pour suppl{\'e}mentaire $Ad(H)$-invariant 
de l'alg{\`e}bre de Lie $\h$ de $\Delta$ dans $\g \oplus 
\g \oplus \g$, le sous-espace vectoriel 
$\m$ form{\'e} des vecteurs {\`a} composante tangente au {\it premier} 
facteur nulle de telle sorte que la restriction de $g \times g \times
g$ {\`a} $\m$ donne une m{\'e}trique sur $M$ qui n'est pas la m{\'e}trique produit
de $G \x G$ mais le rend naturellement r{\'e}ductif. Maintenant la permutation 
circulaire de $G \times G \times G$ est un automorphisme d'ordre 3 qui induit 
une structure d'espace 3-sym{\'e}trique sur $M$ et, finalement, une structure NK.

Lorsque $G=SU(2) \simeq S^3$, si on joint deux bases
orthonorm{\'e}es de $\mathfrak{g}$ et qu'on prolonge par invariance {\`a}
gauche, on obtient un rep{\`e}re 
sur $M$ non orthonorm{\'e} mais dans le co-rep{\`e}re associ{\'e} 
$(e_1,e_2,e_3,f_1,f_2,f_3)$ duquel la forme de K{\"a}hler 
s'{\'e}crit
\begin{equation}
\omega = \frac{\sqrt 3}{2}(e_1 \wedge f_1 + e_2 \wedge f_2 + e_3
\wedge f_3)
\label{can}
\end{equation}

Notre but est de d{\'e}montrer que la structure SNK ainsi
d{\'e}crite sur $S^3 \times S^3$ est la seule (bien s{\^u}r cela d{\'e}pend de
la m{\'e}trique biinvariante choisie, mais comme $SU(2)$ est simple
elles sont toutes proportionnelles):

\begin{prop}
Soit $(g,J)$ une structure presque hermitienne invariante telle que $(S^3 \x
S^3,g,J)$ est SNK. Elle est isomorphe, en tant qu'espace homog{\`e}ne
presque hermitien, {\`a} la construction de Ledger et Obata. En
particulier elle est 3-sym{\'e}trique.
\label{S3}
\end{prop}

\begin{defi}
{\'E}tant donn{\'e} un co-rep{\`e}re $(e_1,e_2,e_3,f_1,f_2,f_3)$ de $S^3
\times S^3$, on appelle 2-forme canonique, la forme donn{\'e}e en
(\ref{can}). 
\end{defi}

Les vari{\'e}t{\'e}s SNK de dimension 6 sont
caracteris{\'e}es par la v{\'e}rification par la forme de K{\"a}hler 
d'une {\'e}quation diff{\'e}rentielle donn{\'e}e par Reyes Carri{\'o}n au th{\'e}or{\`e}me
4.9 page 48 de \cite{re}
(pour une pr{\'e}sentation et une demonstration diff{\'e}rentes du m{\^e}me
fait voir aussi \cite{hi}.) 
\begin{equation}
-2\lambda\omega \wedge \omega = d \hat \rho
\label{equadiff}
\end{equation}
o{\`u} $\lambda$ est une constante r{\'e}elle, $\rho$ est proportionnelle 
{\`a} la diff{\'e}rentielle de $\omega$ et quelle que soit la 3-forme
$\alpha$ de type $(3,0)+(0,3)$, $\hat \alpha$ est l'unique 3-forme
telle que $\alpha + i\hat \alpha$ est une 3-forme volume
complexe. Plus pr{\'e}cisement on notera
\begin{equation}
d\omega=3\lambda \rho
\label{rho}
\end{equation}
D{\`e}s lors on cherche des 2-formes invariantes sur $S^3 \times
S^3$ v{\'e}rifiant (\ref{equadiff}) ou encore des triplets
$(\omega, \rho, \lambda)$ v{\'e}rifiant (\ref{equadiff}) et (\ref{rho}).

\vs

Pour faciliter les calculs on introduit un rep{\`e}re appropri{\'e} dans
lequel l'expression de la forme de K{\"a}hler est proche de celle de la
2-forme canonique.

Sur la sph{\`e}re $S^3$ il existe un rep{\`e}re global $(X_1,X_2,X_3)$
privil{\'e}gi{\'e}, invariant {\`a} gauche et v{\'e}rifiant
\begin{eqnarray*}
   [X_1,X_2] & = & 2X_3 \\
\  [X_2,X_3] & = & 2X_1 \\
\  [X_3,X_1] & = & 2X_2 
\end{eqnarray*}
En prenant la base duale de l'espace cotangent en chaque point
on obtient un rep{\`e}re des 1-formes $(e_1,e_2,e_3)$ tel que
\begin{equation}
de_i=e_{i+1} \wedge e_{i+2} \label{dei}
\end{equation}
en notant les indices dans $\ZM /3 \ZM$. On remarque que n'importe quel
co-rep{\`e}re obtenu a partir d'un tel rep{\`e}re par une isom{\'e}trie directe a
encore la propri{\'e}t{\'e} (\ref{dei}). 

\begin{defi} 
On appelle co-rep{\`e}re circulaire un co-rep{\`e}re $(e_1,e_2,e_3,f_1,f_2,f_3)$
de $S^3 \x S^3$ tel que les trois premi{\`e}res (resp. les trois
derni{\`e}res) formes sont nulles en chaque point sur l'espace
tangent a la fibre du second (resp. du premier) facteur passant par ce
point et v{\'e}rifient (\ref{dei}).
\end{defi}

C'est dans un tel rep{\`e}re que les calculs seront faits. 

\begin{lemm}
Dans un co-rep{\`e}re circulaire $(e_1,e_2,e_3,f_1,f_2,f_3)$, la forme de
K{\"a}lher d'une vari{\'e}te NK s'{\'e}crit seulement comme combinaison
lin{\'e}aire de termes \og mixtes \fg \, $e_i \wedge f_j \ i,j=1,2,3$
\label{C}
\end{lemm}
 
\begin{proof}
On {\'e}crit pour l'instant en toute g{\'e}n{\'e}ralit{\'e}
\[ \omega = \sum_{i=1}^3 a_i e_{i+1} \wedge e_{i+2} + \sum_{i=1}^3 b_i
f_{i+1} \wedge f_{i+2} + \sum_{i,j=1}^3 c_{i,j}
e_i \wedge f_j \]
Que les coefficients soient des constantes, non des fonctions sur la
vari{\'e}t{\'e}, traduit l'invariance de $\omega$.
Puis, en exprimant seulement que la forme est non d{\'e}g{\'e}ner{\'e}e, soit
$\omega \wedge \omega \wedge \omega \neq 0$, on obtient d'abord
\[ ^tACB + det\ C \neq 0
\]
o{\`u} $A$ est le vecteur colonne des $a_i,\ i=1,2,3$, $B$ le
vecteur colonne des $b_i$ et $C$ la matrice des $c_{i,j}$.
Mais pour une vari{\'e}t{\'e} NK, $d\omega$ est une
3-forme de type $(3,0) + (0,3)$. Comme $\omega$ est elle-m{\^e}me de type
$(1,1)$, on ne peut qu'avoir que
\begin{equation}
\omega \wedge d\omega = 0
\end{equation}
c'est {\`a} dire que $\omega \wedge \omega$ est ferm{\'e}e. En fait elle
est m{\^e}me exacte par (\ref{equadiff}). Cela conduit {\`a}
\[ ^tCA = 0 \quad \text{et} \quad CB = 0 \]
puis {\`a}
\[ ^tACB = 0 \]
Finalement $det \ C \neq 0$, c'est {\`a} dire $C$ et $^tC$ sont inversibles et il faut
$A=B=0$.
\end{proof}

Remarquons que l'{\'e}nonc{\'e} est toujours valide si on remplace ``NK'' par ``semi-k{\"a}hlerienne''. 

On peut m{\^e}me mieux choisir sa base pour que $C$ soit diagonale : 
\begin{lemm}
Il existe un co-rep{\`e}re circulaire $(e_1,e_2,e_3,f_1,f_2,f_3)$ tel que
\begin{equation}
\omega = \lambda_1 \, e_1 \wedge f_1 + \lambda_2 \, e_2 \wedge f_2 + 
\lambda_3 \, e_3 \wedge f_3
\label{diag}
\end{equation}
o{\`u} les $\lambda_i \ i=1,2,3$ sont des constantes r{\'e}elles non nulles.
\end{lemm}
\begin{proof}
Soient $M$, $N$ deux matrices de $SO(3)$. 
Partant d'un co-rep{\`e}re circulaire $(e_1,e_2,e_3,f_1,f_2,f_3)$ on lui
associe le co-rep{\`e}re obtenu en appliquant sur chaque espace tangent
aux 3 premi{\`e}res formes la matrice de changement de base $M$ et aux 3
derni{\`e}res la matrice $N$. Comme remarqu{\'e}
pr{\'e}c{\'e}demment, c'est encore un co-rep{\`e}re circulaire (cela
revient {\`a} appliquer sur chaque facteur du produit une isom{\'e}trie directe de la
m{\'e}trique biinvariante de $SU(2)$.) Alors si $\omega$ s'{\'e}crivait dans
l'ancien co-rep{\`e}re gr{\^a}ce {\`a} la matrice $C$, son expression dans le
nouveau fait intervenir la matrice $MC^tN$. En {\'e}crivant
$C$ comme produit d'une matrice sym{\'e}trique (donc diagonalisable par un
changement de base orthonorm{\'e}e) et d'une matrice orthogonale on voit
qu'on peut choisir $M$ et $N$ telles qu'elle soit diagonale.
\end{proof}

\vs

A partir seulement de la forme de K{\"a}hler d'une vari{\'e}t{\'e} SNK, on
peut retrouver la structure presque complexe, puis la
m{\'e}trique.
\begin{prop}
Soit une vari{\'e}t{\'e} presque hermitienne de dimension 6, $(M,g,J)$ munie d'une d'une 3-forme
$\rho$ de type $(3,0)+(0,3)$. Il existe une constante positive $c$ telle que 
\begin{equation}
\iota (v) \rho \wedge \rho=c\ \iota(Jv)vol
\label{J}
\end{equation}
o{\`u} $vol$ d{\'e}signe la forme volume de $g$.
\end{prop}

\begin{proof}
Soit $(v_0,v_1,v_2,v_3,v_4,v_5)$ un rep{\`e}re de $M$ adapt{\'e} {\`a} sa 
structure presque hermitienne : orthonorm{\'e} et tel que
$Jv_{2i}=v_{2i+1}$ pour tout $i=0,1,2$. Autrement dit si
$(l_0,l_1,l_2,l_3,l_4,l_5)$ est le rep{\`e}re \og dual \fg \, 
des 1-formes, la forme de K{\"a}hler est
\[ \omega = l_0 \wedge l_1 + l_2 \wedge l_3 + l_4 \wedge l_5
\]
La forme volume est $vol=l_0 \wedge l_1 \wedge l_2 \wedge l_3 
\wedge l_4 \wedge l_5$. Maintenant $\Lambda^{3,0}$ est 
de dimension complexe 1 donc il existe $z = a 
+ ib \in \CM$ tel que $\rho$ est la partie r{\'e}elle de 
$z(l_0+il_1) \wedge (l_2+il_3) \wedge (l_4+il_5)$ :
\begin{eqnarray*}
\rho & = & a(l_0 \wedge l_2 \wedge l_4 - l_1 \wedge l_3 \wedge l_4 -
l_0 \wedge l_3 \wedge l_5 - l_1 \wedge l_2 \wedge l_5) \\
& & \qquad - b(l_0 \wedge l_2 \wedge l_5 + l_0 \wedge l_3 \wedge l_4 
- l_1 \wedge l_3 \wedge l_5 + l_1 \wedge l_2 \wedge l_4)
\end{eqnarray*}
Le stabilisateur en chaque point d'une telle 3-forme est $SL(3,
\CM)$. Il agit transitivement sur les vecteurs et
pr{\'e}serve (\ref{J}). Il suffit alors de v{\'e}rifier
\[ \iota(l_0)\rho = a(l_2 \wedge l_4 - l_3 \wedge l_5) - b(l_2 \wedge
l_5 + l_3 \wedge l_4) 
\]
puis 
\[ \iota(l_0)\rho \wedge \rho = 2(a^2 + b^2)\iota(l_1)vol \] 
\end{proof}

Cette expression intrins{\`e}que de $J$ permet de le
calculer dans n'importe quel rep{\`e}re : soit $(X_1,X_2,X_3,Y_1,Y_2,Y_3)$
le rep{\`e}re associ{\'e} {\`a} un co-rep{\`e}re circulaire
$(e_1,e_2,e_3,f_1,f_2,f_3)$ dans lequel la forme de K{\"a}hler s'{\'e}crit
(\ref{diag}).
Appelons, en reprenant les notations de Hitchin dans l'article
\cite{hi}, $K_{\rho}$ l'application lin{\'e}aire de $TM$ dans 
$TM \otimes \Lambda^6 T^*M$ d{\'e}finie par
\[ K_{\rho}(v)= \iota (v) \rho \wedge \rho \]
en identifiant $\Lambda^5 T^*M$ et $TM \otimes \Lambda^6 T^*M$, par
l'interm{\'e}diaire du produit ext{\'e}rieur.
On diff{\'e}rentie l'expression (\ref{diag}) de $\omega$ gr{\^a}ce aux 
formules (\ref{dei}). On trouve 
\begin{equation}
K_{d\omega}(X_1)=\big((\lambda_1^2-\lambda_2^2-\lambda_3^2)X_1
-2\lambda_2\lambda_3 Y_1 \big) \otimes s 
\label{K}
\end{equation}
o{\`u} $s=e_1 \wedge e_2 \wedge e_3 \wedge f_1 \wedge f_2 \wedge f_3$. 
Pour calculer la constante $c$ on se sert du fait que $J^2 = -Id$. En
appliquant $K_{d\omega}$ une nouvelle fois {\`a} la partie vectorielle de
l'expression (\ref{K}) ci-dessus, on trouve
\[ c = \frac{1}{3\lambda}\sqrt{-\lambda_1^4-\lambda_2^4-\lambda_3^4
+2\lambda_1^2\lambda_2^2 + 2\lambda_2^2\lambda_3^2 +2\lambda_1^2\lambda_3^2}
\]
Dans cette base, $J$ admet donc la matrice
\[
\left( \begin{array}{cc}
\ D & \ E \\
-E & -D \\
\end{array} \right)
\]
avec
\[ D = \frac{1}{c}
\left( \begin{array}{ccc}
\lambda_1^2-\lambda_2^2-\lambda_3^2 & 0 & 0 \\
0 & \lambda_2^2-\lambda_3^2-\lambda_1^2 & 0 \\
0 & 0 & \lambda_3^2-\lambda_1^2-\lambda_2^2
\end{array} \right)
\]
et
\[ E = \frac{1}{c}
\left( \begin{array}{ccc}
-2\lambda_2\lambda_3 & 0 & 0 \\
0 & -2\lambda_3\lambda_1 & 0 \\
0 & 0 & -2\lambda_1\lambda_2
\end{array} \right)
\]

\vs

On progresse : il est maintenant possible de calculer $\hat \rho$ 
sans passer par la m{\'e}trique, ce qui va nous permettre de r{\'e}soudre
l'{\'e}quation diff{\'e}rentielle (\ref{equadiff}). Quels que soient
les champs de vecteurs $X,Y,Z$ sur $M$ on a 
\[ \hat \rho (X,Y,Z)= -\rho (JX,JY,JZ)= \rho (JX,Y,Z)
\]
\begin{lemm}
Soit $(S^3 \times S^3,\omega,J)$ une vari{\'e}t{\'e} NK et
soit $(e_1,e_2,e_3,f_1,f_2,f_3)$ un co-rep{\`e}re circulaire tel que
$\omega$ s'{\'e}crit (\ref{diag}). Alors (\ref{equadiff}) est
{\'e}quivalente {\`a} 
\begin{equation}
k\frac{1}{\lambda_i}=\lambda_i(\lambda_i^2-\lambda_{i+1}^2-\lambda_{i+2}^2)
\quad \forall i=1,2,3
\label{li2}
\end{equation}
avec 
\begin{equation}
k=\frac{6 \lambda^2 det \ C}{c}
\end{equation} 
\end{lemm}

\begin{proof}
Soit $(e_1,e_2,e_3,f_1,f_2,f_3)$ un co-rep{\`e}re circulaire dans lequel
$\omega$ s'{\'e}crit, comme d{\'e}montr{\'e} au lemme \ref{C},
$\sum c_{i,j} \, e_i \wedge f_j$ et tel que $J$ est
repr{\'e}sent{\'e} dans le rep{\`e}re associ{\'e} par la matrice
\[
\left( \begin{array}{cc}
\ D & \ E \\
F & G \\
\end{array} \right)
\]
On a
\[ \hat{d\omega} (e_{i+1},e_{i+2},f_j)= (C^tD)_{i,j} \]
\[ \hat{d\omega} (e_i,f_{j+1},f_{j+2})= (GC)_{i,j} \]
puis
\[ 3 \lambda \, d\hat \rho = \sum_{i,j=1}^3 (C^tD + GC)_{i,j} \ e_{i+1}
\wedge e_{i+2} \wedge f_{i+1} \wedge f_{i+2} 
\]
Quant au calcul de $\omega \wedge \omega$ il fait appara{\^\i}tre les
mineurs d'ordre 2 de $C$, c'est {\`a} dire son inverse. L'{\'e}quation 
caract{\'e}ristique (\ref{equadiff}) s'{\'e}crit alors
\[ 2k'\,^tC^{-1}=C^tD + GC \]
avec
\begin{equation}
k'=6 \lambda^2 det \ C
\end{equation}
Ce n'est rien d'autre que ($\ref{li2}$), en tenant compte de nos
simplifications successives.
\end{proof}

On r{\'e}soud facilement (\ref{li2}). Premi{\`e}rement notons que si
tous les $\lambda_i$ sont {\'e}gaux ou m{\^e}me seulement de
signes diff{\'e}rents les trois {\'e}quations sont v{\'e}rifi{\'e}es {\`a} la fois pour
$k=\lambda_1^4=\lambda_2^4=\lambda_3^4$.
Autrement si on note $S=\lambda_1^2 + \lambda_2^2 + \lambda_3^2$ \, 
les $\lambda_i^2$ doivent d'abord tous v{\'e}rifier la m{\^e}me {\'e}quation 
du second ordre
\[ 2x^2-Sx-k' =0
\]
Puis supposons que $\lambda_1^2$ et $\lambda_2^2$ soient deux racines
distinctes et, par exemple,
$\lambda_3^2=\lambda_2^2$. Alors $\lambda_1^2\lambda_2^2$ vaut le produit 
des racines : $-\frac{k}{2}$ et l'{\'e}quation (\ref{li2}), $i=2$ implique 
$k=0$, ce qui signifierait que $\omega$ est d{\'e}g{\'e}n{\'e}r{\'e}e.
Par cons{\'e}quent toutes les valeurs diagonales sont {\'e}gales au signe pr{\`e}s.

\vs

Il reste {\`a} trancher cette ambigu{\"\i}t{\'e}. On le fait en introduisant
pour la premi{\`e}re fois la m{\'e}trique, en demandant qu'elle soit positive.

Remarquons premi{\`e}rement que les cas o{\`u} les trois signes sont
positifs ou o{\`u} seulement un signe sur trois l'est sont
identiques, {\`a} une rotation d'angle $\pi$ pr{\`e}s. De m{\^e}me les
deux cas restant. 

Effectivement, si on {\'e}tudie la forme quadratique $X \mapsto
\omega(JX,X)$, on voit qu'elle est soit d{\'e}finie positive si le
d{\'e}terminant de $C$ est positif, soit d{\'e}finie n{\'e}gative dans
le cas contraire. 

\begin{prop}
Soit $\omega$ une 2-forme diff{\'e}rentielle sur $S^3 \times S^3$
v{\'e}rifiant l'{\'e}quation diff{\'e}rentielle (\ref{equadiff}). Il existe
un co-rep{\`e}re circulaire tel que $\omega$ est un multiple de la 2-forme 
canonique. Il est strictement positif si et seulement si $\omega$
repr{\'e}sente une vari{\'e}t{\'e} riemannienne NK.
\end{prop}

\begin{proof}
Pour tout $i$ dans $\ZM /3\ZM$
\[
\omega(Je_i,e_1)=-\frac{2}{c}\lambda_{i+1}\lambda_{i+2}
\, \omega(f_1,e_1)
=\frac{2}{c}\lambda_1\lambda_2\lambda_3
\]
\[
\omega(Je_i,f_i)=\frac{1}{c}\lambda_i(\lambda_i^2-\lambda_{i+1}^2
-\lambda_{i+2}^2)
\]
La forme quadratique $X \mapsto \omega(JX,X)$ est la somme de trois formes
quadratiques de degr{\'e} 2, 
\[ q_i = \frac{2\lambda_i}{c}\big(\lambda_{i+1}\lambda_{i+2}x_i^2-
(\lambda_i^2-\lambda_{i+1}^2-\lambda_{i+2}^2)x_i y_i + 
\lambda_{i+1}\lambda_{i+2}y_i^2\big), \quad i=1,2,3
\]
dont le discriminant est $c^2$, positif, et les coefficients des
termes carr{\'e}s ont le signe de $det\ C$.
\end{proof}

Ceci ach{\`e}ve en m{\^e}me temps la preuve de la proposition \ref{S3}.

\section{Espaces homog{\`e}nes quotients de groupes produits des sph{\`e}res}

Soit $(\g,\h)$ un couple d'alg{\`e}bres de Lie admis dans la liste du
lemme \ref{liste}. On cherche les plongements de $\h$
dans $\g$. En travaillant avec les alg{\`e}bres de Lie on {\'e}carte
provisoirement la question des quotients finis des groupes. Cependant
si $G$ est de la forme $(S^1)^p \x G'/\Gamma$ et $H = (S^1)^q \x  H' /
\Sigma$ o{\`u} $\Gamma$ et $\Sigma$ sont des groupes finis et $G'$ et $H'$ sont
simplement connexes, on retient que la
surjectivit{\'e} de $\phi$ vue en (\ref{istar}) implique premi{\`e}rement
$p \leq q$, deuxi{\`e}mement que le morphisme de groupe obtenu en restreignant
au facteur $(S^1)^p$ de $H$ et en projetant dans $G$ sur le facteur
$(S^1)^q$ est lui-m{\^e}me surjectif. Au niveau des alg{\`e}bres de Lie cela
se traduit, en notant $p$ la projection sur $\bigoplus i\RM$,
parall{\`e}lement {\`a} $\bigoplus \su(2)$, par 
\begin{equation}
\h \stackrel{\varphi}{\hookrightarrow} \g \stackrel{p}{\to} \bigoplus
i\RM \quad \text{est surjectif}
\label{iop}
\end{equation}
 
Dans cette section on s'int{\'e}resse aux espaces homog{\`e}nes $G/H$ o{\`u} $G$
et $H$ sont des produits directs de $S^1$ et $S^3$ ou des
quotients finis de ces produits. Par cons{\'e}quent $\h$ doit {\^e}tre
$i\RM$, $i\RM \oplus i\RM$, $\su(2)$ ou $i\RM \oplus \su(2)$ et d'apr{\`e}s
le lemme \ref{liste} on a toujours $\g = \h \oplus \su(2) \oplus
\su(2)$. Cependant on ignore encore si $\su(2) \oplus \su(2)$
aparaissant dans cette d{\'e}composition est $Ad(H)$-invariant.

\begin{prop}
Soit $M=G/H$ un espace homog{\`e}ne NK simplement connexe de dimension 6 avec $\g$ et $\h$ sommes directes d'alg{\`e}bres de Lie isomorphes {\`a} $i\RM$ ou $\su(2)$. Alors
$\h$ admet pour suppl{\'e}mentaire dans $\g$ un id{\'e}al $\m$
isomorphe {\`a} $\su(2) \oplus \su(2)$ et $M$ est isomorphe {\`a} $S^3
\times S^3$ muni de son unique structure NK 3-sym{\'e}trique invariante. 
\label{iR+su2}
\end{prop}

\begin{proof}
Un id{\'e}al est bien s{\^u}r en particulier $ad(\h)$-invariant et m{\^e}me
$Ad(H)$-invariant puisque $H$ est connexe.
Consid{\'e}rons chaque cas :

\vs
-- $\g=i\RM \oplus \su(2) \oplus \su(2), \h = i\RM$

D'apr{\`e}s (\ref{iop}), $\varphi \circ p$ est surjectif, c'est {\`a} dire dans ce
cas bijectif. L'intersection du plongement de $\h$ avec le
noyau de la projection est $\varphi(\h) \cap (\su(2) \oplus
\su(2))= \{0\}$ et ce dernier est par cons{\'e}quent toujours un suppl{\'e}mentaire de
$\h$ dans $\g$, quel que soit pr{\'e}cis{\'e}ment le plongement.

\vs
-- $\g=i\RM \oplus i\RM \oplus \su(2) \oplus \su(2), \h = i\RM
\oplus i\RM$

De m{\^e}me ici $\m = \su(2) \oplus \su(2)$ convient.

\vs
-- $\g=\su(2) \oplus \su(2) \oplus \su(2), \h = \su(2)$

En projetant sur chaque facteur $\su(2)$ de $\g$ on obtient un
endomorphisme de $\su(2)$, soit une repr{\'e}sentation unitaire de
dimension 2 de $\su(2)$. Or $\su(2)$ n'a qu'une seule repr{\'e}sentation
irr{\'e}ductible complexe en chaque dimension. En dimension 2 on ne peut
donc avoir que l'identit{\'e} de $\su(2)$ ou la repr{\'e}sentation triviale. Si
les trois repr{\'e}sentations {\'e}taient triviales ce ne serait pas un
plongement. Il y a donc au moins un facteur $\su(2)$ tel que la
projection $q$ v{\'e}rifie comme plus haut : $\varphi \circ q$ est surjectif
et on prend $\m$ {\'e}gal au noyau de $q$ : la somme des deux
autres. C'est encore un id{\'e}al.

\vs
-- $\g=i\RM \oplus \su(2) \oplus \su(2) \oplus \su(2), \h = i\RM \oplus \su(2)$

Cela r{\'e}sulte de la combinaison des deux arguments pr{\'e}c{\'e}dents
pour les facteurs $i\RM$ et $\su(2)$, respectivement, de $\g$.

\vs
Un id{\'e}al est m{\^e}me une sous-alg{\`e}bre de Lie. La vari{\'e}t{\'e} {\'e}tant
simplement connexe est isomorphe {\`a} $S^3 \x S^3$. De plus puisque $\m
\simeq \su(2) \oplus \su(2)$, on peut toujours le munir de
l'unique ({\`a} multiple pr{\`e}s) 2-forme $\omega$ susceptible de repr{\'e}senter
une structure NK invariante d'apr{\`e}s les r{\'e}sultats de la section
3. Mais ici $H$ n'est plus trivial et il faut que $\omega$ soit
$Ad(H)$-invariante pour le plongement pr{\'e}cis de $H$ choisi. Par
construction, c'est le cas d{\`e}s que $Ad(H)$ est inclus dans le sous-groupe
diagonal. C'est de plus une condition n{\'e}c{\'e}ssaire, $Ad(H)$ agissant
s{\'e}par{\'e}ment sur chaque facteur de la somme. Les seules possibilit{\'e}s
pour $H$ sont finalement les sous-groupes de $S^3$ : $\{1\}$, $S^1$ et
les quotients finis de $S^3$. En effet
un espace homog{\`e}ne qu'on {\'e}crit $M=G/H$ peut toujours s'{\'e}crire
diff{\'e}remment $M=G'/H'$ o{\`u} $G'$ est un sous-groupe 
d'isom{\'e}tries plus petit que $G$ mais agissant toujours
transitivement et $H'$ le sous-groupe d'isotropie {\it dans} $G'$.

\end{proof}

\section{L'espace des drapeaux}

L'espace des drapeaux $F(1,2)$ d'un espace vectoriel hermitien $E$ de dimension
$3$ est l'espace des couples $(l,p)$ o{\`u} $l$ est une droite de $E$ et
$p$ un plan contenant cette droite. Il appara{\^\i}t naturellement dans
notre liste comme $SU(3)/S^1 \x S^1$ mais on peut aussi l'{\'e}crire
$U(3)/U(1) \x U(1) \x U(1)$ : un point $(l,p)$ est autrement d{\'e}fini
par une base orthonorm{\'e}e $(e_1, e_2, e_3)$ telle que $l = \CM e_1$ et
$(e_1, e_2)$ est une base orthonorm{\'e}e de $p$ ; l'action naturelle de
$U(3)$ sur $E$ induit une action sur les drapeaux dont le groupe
d'isotropie en un point $(l,p)$ est form{\'e} d'endomorphismes qui
pr{\'e}servent les trois droites complexes $\CM e_1, \CM e_2, \CM e_3$,
c'est {\`a} dire d'endorphismes diagonaux dans la base $(e_1, e_2,
e_3)$. De fa\c con {\'e}quivalente il existe trois fibrations $F(1,2) \to
\CM P(2)$ {\`a} fibres isom{\'e}triques {\`a} $\CM P(1)$. Sur la fibre de la
premi{\`e}re c'est la droite qui varie dans le plan, sur celle de la
seconde c'est le plan autour de la droite et sur la fibre de la
troisi{\`e}me, la deuxi{\`e}me droite du plan, $\CM e_2$, est
fixe et la droite $l$ varie dans le plan orthogonal et le plan $p$
avec elle autour de $e_2$. En fait chacune de ces fibrations est la
fibration d'un espace de twisteurs au-dessus d'une vari{\'e}t{\'e} de
dimension 4. Cela permet que $F(1,2)$ soit muni naturellement de trois
structures k{\"a}hleriennes puis, par variation canonique de la
submersion riemannienne, d'une structure NK (voir \cite{na2} pour la
construction d'une vari{\'e}t{\'e} NK {\`a} partir d'une submersion
k{\"a}hlerienne g{\'e}n{\'e}rale, \cite{bau} pour
les espaces de twisteurs NK.) Pour chercher toutes les
structures NK de $F(1,2)$ on ne privil{\'e}gie aucune fibration ou
aucune direction complexe associ{\'e}e {\`a} un point de l'espace de
drapeaux, on regarde la repr{\'e}sentation complexe $Ad(H)$ de
$H=U(1) \x U(1) \x U(1)$ associ{\'e}e au plongement naturel dans $G=U(3)$.
Tous les plongements de $H$ dans $G$ sont conjugu{\'e}s et induisent le
m{\^e}me espace homog{\`e}ne $G/H$ en fin de compte car leur image
est un tore maximal. On choisit un suppl{\'e}mentaire $\m$, $Ad(H)$ invariant. 
\begin{lemm}
L'espace des m{\'e}triques invariantes presque hermitiennes de $F(1,2)$
est de dimension 3
\end{lemm}
\begin{proof}
La repr{\'e}sentation complexe $Ad(H)$ sur $\m$ (ou la repr{\'e}sentation
lin{\'e}aire isotropique) est r{\'e}ductible : elle se d{\'e}compose en une
somme de trois repr{\'e}sentations irr{\'e}ductibles. En effet si on
repr{\'e}sente habituellement $U(3)$ par les matrices unitaires et le
sous-groupe $H=S^1 \x S^1 \x S^1$ par les matrices diagonales, $\u(3)$
est l'ensemble des matrices anti-hermitiennes et l'ensemble $\k$ des
matrices avec des z{\'e}ros sur la diagonale est un suppl{\'e}mentaire
{\'e}vident de $\h$. Or il est $Ad(H)$-invariant. En fait si on note
\[
\forall \, a,b,c \in \CM \quad \la a,b,c \ra \stackrel{def}{=} 
\left(\begin{array}{ccc}
0             & a             & b \\
-\overline{a} & 0             & c \\
-\overline{b} & -\overline{c} & 0
\end{array}\right)
\]
\begin{equation}
Ad_h \la a,b,c \ra = \la e^{i(t-s)}a,e^{i(t-r)}b,e^{i(s-r)}c \ra \quad \text{o{\`u} } h
= 
\left(\begin{array}{ccc}
e^{ir} & 0 & 0 \\
0 & e^{is} & 0 \\
0 & 0 & e^{it}
\end{array}\right)
\in H
\label{AdH}
\end{equation}
On le scinde en trois sous-espaces invariants :
\begin{eqnarray*}
                \l & = & \{\la a,0,0 \ra \mid a \in \CM \} \\
                \m & = & \{\la 0,b,0 \ra \mid b \in \CM \} \\
\text{et} \quad \n & = & \{\la 0,0,c \ra \mid c \in \CM \}
\end{eqnarray*}
On voit par (\ref{AdH}) que
\begin{equation} 
[\h,\l] = \l, \ [\h,\m] = \m, \ [\h,\n] = \n
\label{crochets}
\end{equation}
La repr{\'e}sentation $Ad(H)$ restreinte {\`a} $\l$, $\m$ ou $\n$ est
irr{\'e}ductible et l'espace des produits scalaires
$Ad(H)$-invariants de $\k$ ou de fa\c con {\'e}quivalente l'espace des m{\'e}triques
invariantes de $F(1,2)$ est de dimension 3.
\end{proof}

De plus on peut calculer
\begin{eqnarray*}
      [ \la a,0,0 \ra,\la 0,b,0 \ra ]       & = & \la 0,0,-\overline{a}b \ra \\
\lbrack \la a,0,0 \ra,\la 0,0,c \ra \rbrack & = & \la 0,ac,0 \ra \\
\lbrack \la 0,b,0 \ra,\la 0,0,c \ra \rbrack & = & \la -b\overline{c},0,0 \ra
\end{eqnarray*}
et
\[ 
[\la a,0,0 \ra,\la a',0,0 \ra] = 
\left(\begin{array}{ccc}
iz & 0   & 0 \\
0  & -iz & 0 \\
0  & 0   & 0
\end{array}\right)
\quad \text{o{\`u} } z=2Im(\overline{a}a')
\]
et de m{\^e}me pour $\m$ et $\n$ si bien que
\begin{eqnarray}
[\l,\l] \subset \h, \  [\m,\m] \subset \h, \ [\n,\n] \subset \h
\label{crochets2} \\
\lbrack\l,\m\rbrack = \n, \ [\m,\n] = \l, \ [\n,\l] = \m \label{crochets3}
\end{eqnarray}
Ces relations sont compatibles avec celles calcul{\'e}es par les
auteurs de \cite{bau} p141 o{\`u} notre \og $\l \oplus \m$ \fg \ est
not{\'e} \og $\m$ \fg. On consid{\`e}re comme eux le produit scalaire sur $\u(3)$
donn{\'e} par :
\[ \la X|Y \ra = -\frac{1}{2}Re(tr(XY)) \]
Les trois sous-espaces invariants sont deux {\`a} deux orthogonaux, on
peut donc appliquer une homot{\'e}tie dans chacun et
obtenir encore une m{\'e}trique invariante. En fait on les obtient toutes
de cette fa\c con. On note $g$ la m{\'e}trique invariante associ{\'e}e aux coefficients
d'homot{\'e}tie $r,s,t \in ]0;+\infty[$ c'est {\`a} dire au produit scalaire
\[ (.|.) = r \la .|. \ra |_{\l \x \l} + 
s \la .|. \ra |_{\m \x \m} + t \la .|. \ra |_{\n \x \n}
\]
C'est la m{\'e}trique d'une submersion riemannienne au-dessus de $\CM
P(2)$ muni de la m{\'e}trique standard si et seulement si deux param{\`e}tres,
supposons $r$ et $s$, sont {\'e}gaux {\`a} $1$. Alors $\n$ est
l'espace tangent {\`a} la fibre en l'origine et (\ref{crochets2})
exprime que les fibres sont totalement g{\'e}od{\'e}siques. Puisqu'il
s'agit de la fibration d'un espace de twisteur au-dessus d'une
vari{\'e}t{\'e} d'Einstein auto-duale de courbure scalaire {\'e}gale {\`a} $24$,
le th{\'e}or{\`e}me de Hitchin (\cite{hi}) ou Friedrich (\cite{fr}) dit que
$(F(1,2),g)$ est k{\"a}hlerienne si et seulement si le troisi{\`e}me param{\`e}tre
$t$ vaut 2. De plus par un autre th{\'e}or{\`e}me de Friedrich dans \cite{fr2}
\begin{theo}
Soit $(M^4,g)$ une vari{\'e}t{\'e} riemannienne de dimension 4. Si son
espace de twisteurs $Z \stackrel{\pi}{\to} M$, muni de la m{\'e}trique
$\pi^*g + tds^2$, o{\`u} $ds^2$ est la m{\'e}trique standard de $\CM
P(1)$, est d'Einstein pour un certain $t>0$, alors $(M^4,g)$ est
auto-duale, d'Einstein, {\`a} courbure scalaire strictement positive $R$,
et $t$ vaut $\frac{48}{R}$ ou $\frac{24}{R}$.
\end{theo}
\noindent elle est d'Einstein si et seulement si $t=1$ ou $t=2$. Dans
le premier cas elle est donc K{\"a}hler-Einstein, dans
le second il s'av{\`e}re (voir \cite{bau}, p145) qu'elle est strictement NK. Comme
en dimension 6 une vari{\'e}t{\'e} NK est k{\"a}hlerienne ou d'Einstein (\cite{gr}),
les seules possibilit{\'e}s que la vari{\'e}t{\'e} soit
NK dans le cas o{\`u} deux param{\`e}tres sont {\'e}gaux, par
exemple $r=s$, sont finalement $t=2r=2s$ et $t=r=s$.
 
A cause de (\ref{crochets}), (\ref{crochets2}), (\ref{crochets3}), on
cherche une application $\Lambda: \k \to \mathfrak{so}(\k)$,
repr{\'e}sentant la connexion de Levi-Civit{\'a} de $g$, de la forme
\begin{eqnarray*}
\Lambda(X)U & = & \alpha [X,U] \\
\Lambda(U)A & = & \beta [U,A] \\
\Lambda(A)X & = & \gamma [A,X] \\
\Lambda(X)Y & = & \Lambda(U)V = \Lambda(A)B = 0
\end{eqnarray*} 
o{\`u} $X,Y \in \l$, $U,V \in \m$ et $A,B \in \n$. Par (\ref{L}) on a alors
\begin{eqnarray*}
\Lambda(U)X & = & (1-\alpha) [U,X] \\
\Lambda(A)U & = & (1-\beta) [A,U] \\
\Lambda(X)A & = & (1-\gamma) [X,A] \\
\end{eqnarray*} 
En prenant dans (\ref{L'}) $X \in \l, Y \in \m, Z \in \n$ puis en permutant
circulairement on obtient les conditions, d'ailleurs suffisantes
\begin{equation}
\left\{ \begin{array}{lll}
\alpha t & = & (1-\gamma)s \\
\beta r  & = & (1-\alpha)t \\
\gamma s & = & (1-\beta)r
\end{array} \right.
\label{abc}
\end{equation}

\begin{lemm}
Les seules structures presque complexes invariantes de l'espace des
drapeaux, compatibles avec $g$, sont celles repr{\'e}sent{\'e}es par
la multiplication par $\pm i$ de chacun des trois nombres complexes
intervenant dans l'{\'e}criture des matrices de $\k$.
\end{lemm}

\begin{proof}
On note de la m{\^e}me fa\c con $J$ une structure presque complexe
invariante de $F(1,2)$ ou l'endomorphisme de carr{\'e} $-1$ de $\k$ qui la
repr{\'e}sente. Si on note
\[ \forall \, a \in \CM, \quad J \la a,0,0 \ra  = \la a',b',c' \ra \]
\[ Ad_h J \la a,0,0 \ra = \la e^{i(t-s)}a',e^{i(t-r)}b',e^{i(s-r)}c' \ra \quad
\text{pour } h = 
\left(\begin{array}{ccc}
e^{ir} & 0 & 0 \\
0 & e^{is} & 0 \\
0 & 0 & e^{it}
\end{array}\right)
\in H
\]
Dans l'autre sens
\[
Ad_h \la a,0,0 \ra = \la e^{i(t-s)}a,0,0 \ra
\]
On choisit de faire $s=t$.
Alors il faut quels que soient $r$ et $s$
\[
\la a',b',c' \ra = J Ad_h \la a,0,0 \ra = Ad_h J \la a,0,0 \ra 
= \la a',e^{i(s-r)}b',e^{i(s-r)}c' \ra
\]
car l'{\'e}galit{\'e} centrale doit {\^e}tre vraie quel que soit $h \in H$.
La seule solution est que $b'=c'=0$ quel que soit $a \in \CM$.
On proc{\`e}de de la m{\^e}me mani{\`e}re pour les deux autres sous-espaces et on trouve
que $J$ pr{\'e}serve $\l$, $\m$ et $\n$. Sur ces sous-espaces de dimension
2, ce ne peut-{\^e}tre que la rotation d'angle $\pm \frac{\pi}{2}$ par rapport {\`a}
$(.|.)|_{\l \x \l}, (.|.)|_{\m \x \m}$ ou $(.|.)|_{\n \x \n}$. 
\end{proof}
On note dans la suite $\e_1, \e_2, \e_3$ les nombres, {\'e}gaux {\`a} $\pm 1$, tels que
\[
J \la a,b,c \ra = \la \e_1a,-\e_2b,e_3c \ra
\]
Puisqu'on connait la connexion de Levi-Civit{\'a} $\nabla$ on peut
calculer $\nabla J$. En exprimant qu'il doit {\^e}tre antisym{\'e}trique :
\[ (\nabla_X J)Y = -(\nabla_Y J)X \]
on obtient des conditions sur $\e_1, \e_2, \e_3$ et $\alpha, \beta,
\gamma$. Les trois sous-espaces {\'e}tant pr{\'e}serv{\'e}s par $J$,
$(\nabla_X J)Y=0$ d{\`e}s que $X,Y$ appartiennent au m{\^e}me. Les autres
cas donnent :
\begin{equation}
\left\{ \begin{array}{lll}
\alpha (\e_2 + \e_3) & = & (1-\alpha)(\e_1 + \e_3) \\
\beta (\e_1 + \e_3)  & = & (1-\beta)(\e_1 + \e_2) \\
\gamma (\e_1 + \e_2) & = & (1-\alpha)(\e_2 + \e_3) 
\end{array} \right.
\label{+-}
\end{equation}
Si les trois signes sont {\'e}gaux $\alpha = 1-\alpha$,
c'est {\`a} dire $\alpha = \frac{1}{2}$ puis $\beta = \gamma =
\frac{1}{2}$. En reportant dans (\ref{abc}) on trouve $r=s=t$, c'est
{\`a} dire que $g$ est un multiple strictement positif de la m{\'e}trique NK connue.
S'il y a deux signes distincts en revanche toutes les lignes
de (\ref{+-}) sont nulles c'est {\`a} dire $\nabla J = 0$ et selon
desquels il s'agit un param{\`e}tre dans l'{\'e}criture de la m{\'e}trique
est {\'e}gal {\`a} la somme des deux autres. On appelle ces
m{\'e}triques $g_{\lambda,\mu}$ (quand $\alpha$ est nul : les coefficients
de la m{\'e}trique sont alors $s=\lambda$, $t=\mu$ et $r = \lambda + \mu$),
$g'_{\lambda,\mu}$ (lorsque $\beta = 0$) et $g''_{\lambda,\mu}$ ($\gamma=0$),
pour chaque couple de nombres $(\lambda,\mu)$ strictement positifs.
On a d{\'e}montr{\'e} la

\begin{prop}
L'espace des drapeaux peut-{\^e}tre muni d'une seule m{\'e}trique
strictement NK homog{\`e}ne, trois-sym{\'e}trique, {\`a} un changement
d'{\'e}chelle pr{\`e}s, et des seules m{\'e}triques k{\"a}hleriennes
homog{\`e}nes $g_{\lambda,\mu}$, $g'_{\lambda,\mu}$
et $g''_{\lambda,\mu}$, $\lambda,\mu \in ]0;+\infty[$. Ces derni{\`e}res,
lorsque $\lambda=\mu$, correspondent, {\`a} un changement d'{\'e}chelle pr{\`e}s, {\`a}
la m{\'e}trique k{\"a}hlerienne naturelle de l'espace de twisteur de $\CM
P(2)$, la fibration {\'e}tant r{\'e}alis{\'e}e de trois fa\c cons
diff{\'e}rentes {\`a} partir de $F(1,2)$.
\end{prop}

\section{L'espace projectif complexe de dimension 3}

Dans cette section $\h=i\RM \oplus \su(2)$, $\g=\mathfrak{sp}(2)$,
l'alg{\`e}bre de Lie de $Sp(2)$.

\begin{lemm}
Il y a un seul plongement possible, {\`a} conjugaison pr{\`e}s, de $i\RM
\oplus \su(2)$ dans $\mathfrak{sp}(2)$, donn{\'e} par la composition des
plongements naturels $i\RM \hookrightarrow \mathfrak{sp}(1)$ et
$\su(2) \stackrel{\sim}{\longrightarrow} \mathfrak{sp}(1)$ et du plongement
diagonal $\mathfrak{sp}(1) \oplus  \mathfrak{sp}(1) \hookrightarrow
\mathfrak{sp}(2)$.
\end{lemm}

\begin{proof}
On cherche les plongements $j : i\RM \oplus \su(2) \to
\mathfrak{sp}(2)$. On consid{\`e}re la repr{\'e}sentation complexe $\rho$
de dimension 4 induite par la restriction de $j$ {\`a} $\su(2)$. Elle
commute {\`a} tout {\'e}l{\'e}ment $X$ de l'image de $i\RM$. Comme tout
endomorphisme de $\mathfrak{sp}(2)$, on peut voir X comme un
endomorphisme complexe via l'inclusion $\mathfrak{sp}(2) \subset
\su(4)$. Il est alors diagonalisable et son spectre est de la forme
$\{ \lambda i, \mu i, -\lambda i, -\mu i\}$ o{\`u} $\lambda, \mu \in \RM$
car quel que soit le vecteur propre $u$ de $X$, $ju$ est un vecteur
propre pour la valeur propre oppos{\'e}e. Si $\lambda=\mu$, c'est {\`a} dire
$X$ est la multiplication dans $\HM^2$ par un nombre imaginaire pur, $\rho$ est
la somme directe de deux repr{\'e}sentations irr{\'e}ductibles de dimension 2,
n{\'e}c{\'e}ssairement conjugu{\'e}es. Cependant $j$ n'est alors pas un
plongement, son image {\'e}tant isomorphe {\`a} $\su(2)$. Par cons{\'e}quent
$\lambda \neq \mu$ et on doit 
n{\'e}c{\'e}ssairement avoir $\lambda$ ou $\mu = 0$ sinon $\rho$
pr{\'e}serverait les 4 sous-espaces propres de dimension 1, c'est {\`a}
dire serait triviale. Comme la somme des sous-espaces propres de
$\lambda$ et $-\lambda$ est stable par la multiplication par $i,j,k$,
ce cas est exactement celui d{\'e}crit dans l'{\'e}nonc{\'e} du lemme.
\end{proof}
Dans le cas de groupes simples, l'espace homog{\`e}ne obtenu est
$Sp(2)/S^1 \x SU(2)$, isomorphe {\`a} $\CM P(3)$. 

En g{\'e}n{\'e}ral $\CM P(2n-1)$ est isomorphe {\`a} $Sp(n)/S^1 \x
Sp(n-1)$. En effet le groupe $Sp(n)$ agit transitivement sur $\CM ^{2n}$,
identifi{\'e} {\`a} $\HM ^n$, en pr{\'e}servant les droites complexes et le
groupe d'isotropie en $u \in \CM P(2n-1)$ de l'action
induite est constitu{\'e} d'endomorphismes qui pr{\'e}servent non
seulement $u$ mais $ju$ et agissent comme $Sp(n-1)$ sur l'orthogonal,
vu comme $\HM^{n-1}$. Pour
$n=2$, $Sp(2) \simeq Spin(5)$ agit transitivement sur la sph{\`e}re $S^4$ par
l'interm{\'e}diaire de $SO(5)$ et le groupe d'isotropie est l'image
r{\'e}ciproque, par le rev{\^e}tement {\`a} $2$ feuillets $\pi : Spin(5) \to
SO(5)$, de $SO(4)$ i.e. $Spin(4)$. Par cons{\'e}quent $S^4 \simeq
Sp(2)/Spin(4)$ et on d{\'e}finit une fibration $\CM P(3) \to S^4$ par le
plongement naturel $S^1 \x Sp(1) \hookrightarrow Sp(1) \x Sp(1) \simeq
Spin(4)$. C'est en fait la fibration de l'espace de twisteurs de $S^4$
comme expliqu{\'e} dans \cite{re} p45

D{\`e}s lors si on d{\'e}compose l'alg{\`e}bre de Lie de $Sp(2)$ en
$\mathfrak{sp}(5)=\mathfrak{spin}(4) \oplus \m$, $\m$ est identifi{\'e}
{\`a} l'espace tangent {\`a} l'origine de $S^4$. L'espace tangent {\`a}
l'origine de $\CM P(3)$ est lui identifi{\'e} {\`a} $\m \oplus \n$ o{\`u} $\n$
est un suppl{\'e}mentaire de $i\RM$ dans $\su(2)$ si bien que $\g = \h
\oplus \m \oplus \n$, $\mathfrak{spin}(4)$ {\'e}tant isomorphe {\`a} $\su(2)
\oplus \su(2)$.

Explicitement, on note $a \mapsto a^*, \ a \in \HM$ la conjugaison
quaternionique (soit
pour $u,v \in \CM$, $u + jv \mapsto \overline u - jv$). On repr{\'e}sente
habituellement $\g=\mathfrak{sp}(2)$ dans l'espace des matrices
carr{\'e}es $2 \x 2$ {\`a} coefficients quaternioniques et on pose 
\[
\m=\{
\left(\begin{array}{cc}
0   & a \\
a^* & 0
\end{array}\right)
| \ a \in \HM \}
\]
\[
\n=\{
\left(\begin{array}{cc}
a   & 0\\
0 & 0
\end{array}\right)
| \ a = jx + ky, \ x,y \in \RM \}
\]
Leur somme $\k = \m \oplus \n$ est un suppl{\'e}mentaire de $\h$. Bien
plus, chacun de ces sous-espaces est pr{\'e}serv{\'e} par $Ad(H)$. Un
{\'e}l{\'e}ment $h \in H$ repr{\'e}sent{\'e} par 
$\left(\begin{array}{cc}
e^{i\theta} & 0 \\
0           & u
\end{array}\right)$, $u \in Sp(1)$ agit par l'action adjointe comme
\[ 
Ad_h
\left(\begin{array}{cc}
0   & a \\
a^* & 0
\end{array}\right)
=\left(\begin{array}{cc}
0                    & e^{i\theta}ab \\
b^* a^* e^{-i\theta} & 0
\end{array}\right)
\]
sur les {\'e}l{\'e}ments de $\m$
et
\[ 
Ad_h
\left(\begin{array}{cc}
a & 0 \\
0 & 0
\end{array}\right)
=\left(\begin{array}{cc}
e^{2i\theta}a & 0 \\
0             & 0
\end{array}\right)
\]
sur $\n$. Les deux repr{\'e}sentations de $H$ ainsi d{\'e}crites sont
irr{\'e}ductibles, la premi{\`e}re, de dimension 4, a fortiori car les
repr{\'e}sentations de $Sp(1)$ donn{\'e}es par la multiplication {\`a}
droite ou {\`a} gauche dans $\HM$ le sont. Par cons{\'e}quent $Ad(H)\k \subset
\k$ et $\k = \m \oplus \n$ est une d{\'e}composition
irr{\'e}ductible. L'espace des m{\'e}triques homog{\`e}nes est de dimension 2 et si on
s'autorise un changement d'{\'e}chelle elles ne sont plus d{\'e}crites que
par un seul param{\`e}tre, la courbure scalaire de la fibre,
isomorphe {\`a} $\CM P(1)$, de la fibration riemannienne $\CM P(3) \to
S^4$. Autrement dit ce sont les multiples strictement positifs des
m{\'e}triques twistorielles au dessus de $S^4$, munie de sa m{\'e}trique
standard. Plus pr{\'e}cis{\'e}ment, avec le produit scalaire sur $\k$ donn{\'e} par
$\langle X,Y \rangle = -\frac{1}{2}Re(tr(XY))$ soit par
\[ \langle \left(\begin{array}{cc}
0   & a \\
a^* & 0
\end{array}\right),
\left(\begin{array}{cc}
0   & b \\
b^* & 0
\end{array}\right) \rangle
= Re(ab^*), \quad 
\langle \left(\begin{array}{cc}
a & 0 \\
0 & 0
\end{array}\right),
\left(\begin{array}{cc}
b & 0 \\
0 & 0
\end{array}\right) \rangle
= \frac{1}{2}Re(ab^*)
\]
et $\m$ et $\n$ sont orthogonaux, les m{\'e}triques twistorielles sont
les m{\'e}triques invariantes $g_t$ valant $\langle
.,. \rangle |_{\m \x  \m} + t \langle .,. \rangle |_{\n \x
  \n}$ sur $T_o M \simeq \k$. En effet on a vu que $\m$ est identifi{\'e}
{\`a} $T_o(S^4)$. C'est lui qui re\c coit la m{\'e}trique de la base $S^4$. 

Finalement, {\`a} cause des m{\^e}mes
th{\'e}or{\`e}mes de Friedrich et Hitchin sur les espaces de twisteurs en
dimension 6 cit{\'e}s {\`a} la section pr{\'e}c{\'e}dente
\begin{prop}
Dans la famille $(g_t)_{t>0}$ de m{\'e}triques homog{\`e}nes de $\CM P(3)
= Sp(2)/S^1 \x Sp(1)$, seules $g_1$, strictement NK 3-sym{\'e}trique,
et $g_2$, K{\"a}hler-Einstein, sont des m{\'e}triques NK. Toutes les
m{\'e}triques NK homog{\`e}nes sur  $\CM P(3)$ sont des multiples de celles-ci.
\end{prop}

\section{La sph{\`e}re de dimension six}

Il y a un seul espace homog{\`e}ne $M=G_2/SU(3)$ -- c'est {\`a} dire un seul
sous-groupe de $G_2$ {\`a} conjugaison pr{\`e}s isomorphe {\`a} $SU(3)$ -- car
les deux groupes sont de m{\^e}me rang. Il est isomorphe {\`a} la sph{\`e}re de
dimension 6. On peut construire une structure NK 3-sym{\'e}trique sur $M$ en
prenant comme dans \cite{re} sur l'espace tangent l'unique m{\'e}trique,
{\`a} un multiple pr{\`e}s, et l'unique structure presque complexe
pr{\'e}serv{\'e}es par le groupe d'isotropie. En effet $SU(3)$ n'admet pas de
repr{\'e}sentation autre que triviale avant la dimension 6, $G_2/SU(3)$
est donc {\`a} isotropie irr{\'e}ductible. On a par cons{\'e}quent
\begin{prop}
La sph{\`e}re $S^6$, vue comme l'espace homog{\`e}ne $G_2/SU(3)$, admet une
seule structure presque hermitienne invariante {\`a} homot{\'e}tie
pr{\`e}s. Elle est NK, 3-sym{\'e}trique.
\end{prop}

En revanche la vari{\'e}t{\'e} $S^6$, munie de sa m{\'e}trique ronde $g$,
a une infinit{\'e} de structures presque complexes NK. Elles sont
toutes invariantes sous l'action d'un sous-groupe diff{\'e}rent, isomorphe
{\`a} $G_2$, du groupe d'isom{\'e}tries $SO(7)$ c'est {\`a} dire correspondent
chacune {\`a} une fa\c con diff{\'e}rente de r{\'e}aliser $S^6$
comme l'espace homog{\`e}ne presque hermitien $G_2/SU(3)$. En effet
contrairement {\`a} ce qui se passe
pour les autres vari{\'e}t{\'e}s NK comme $F(1,2)$, $\CM P(3)$, $S^3
\x S^3$, il n'y a pas ici une seule droite de spineurs de Killing qui
permettrait de d{\'e}finir l'unique structure
presque complexe compatible avec $g$ rendant la vari{\'e}t{\'e} NK (voir
proposition \ref{spineurs}).
Cependant, soit une vari{\'e}t{\'e} riemannienne $(M,g)$, C. B{\"a}r a montr{\'e},
dans le cadre d'une explication g{\'e}n{\'e}rale \cite{ba} {\`a} la fois de la
th{\'e}orie des spineurs de Killing et de
l'holonomie sp{\'e}ciale, que la donn{\'e}e d'une
structure presque complexe $J$ sur $M$ telle que $(M,g,J)$ soit NK est
{\'e}quivalente {\`a} la donn{\'e}e sur son c{\^o}ne
riemannien d'une 3-forme g{\'e}n{\'e}rique parall{\`e}le $\alpha$ (i.e. il est {\`a} holonomie
contenue dans $G_2$). La correspondance est donn{\'e}e par
\begin{equation}
\alpha = t^2 dt \wedge \omega + t^3 d\omega
\label{3forme}
\end{equation}
o{\`u} $\omega$ est la forme de K{\"a}hler. Ici le c{\^o}ne est
simplement l'espace plat de dimension 7 et la forme parall{\`e}le
est constante. Il y a donc une injection de l'espace des structures
presque complexes NK de $(S^6,g)$ dans l'espace des 3-formes
g{\'e}n{\'e}riques en dimension 7, ouvert dans
$\Lambda^3(\RM^7)$ (voir \cite{br}). 

La d{\'e}riv{\'e}e covariante pour la connexion de Levi-Civit{\'a} et la
diff{\'e}rentielle ext{\'e}rieure de la forme de K{\"a}hler d'une vari{\'e}t{\'e} NK sont
proportionnelles : $d \omega = 3 \nabla \omega$. Les vari{\'e}t{\'e}s SNK de
dimension 6 v{\'e}rifient en outre, pour tous champs de vecteurs $X,Y$ :
\begin{equation}
\| (\nabla_X J)Y \|^2 = \alpha \left( \|X\|^2 \|Y\|^2 - g(X,Y)^2 -
g(JX,Y)^2 \right)
\label{alpha} 
\end{equation}
o{\`u} $\alpha$ est un nombre r{\'e}el strictement positif reli{\'e} {\`a} la
courbure scalaire par $s=30 \alpha$. La vari{\'e}t{\'e} NK est alors dite
\og de type constant $\alpha$ \fg \, (cf \cite{gr2}). Si on fixe la
m{\'e}trique standard sur la sph{\`e}re, de courbure scalaire $s=30$, et la
norme standard de l'espace des 3-formes associ{\'e}e {\`a} cette
m{\'e}trique, $\alpha$ vaut $1$ quelle que soit la structure presque
complexe NK et la norme de $d\omega$ par cons{\'e}quent vaut toujours $6$
par (\ref{alpha}).
   
Alors, on se donne un endomorphisme $J$ de carr{\'e} $-1$ d'un espace tangent
en $x \in S^6$ tel que $g_x(JX,JY)=g_x(X,Y)$, ou de fa\c con {\'e}quivalente
la 2-forme $\omega_x$, et une 3-forme $\rho$ appartenant {\`a} la
sph{\`e}re de dimension 1, de rayon $6$ de
l'espace des 3-formes de type (3,0)+(0,3) de $T_xS^6$. La 3-forme
constante $\alpha$ sur $\RM^7$ donn{\'e}e par $\alpha_x = dt \wedge
\omega_x + \rho$ permet une r{\'e}duction {\`a} $G_2$ de l'holonomie du c{\^o}ne,
ce qui prouve r{\'e}ciproquement que
\[ \omega = \frac{1}{t^2}\iota(dt)\alpha \]
est la 2-forme de K{\"a}hler d'une structure NK sur $S^6$ avec $\rho =
(d\omega)_x$.

\begin{prop}
L'ensemble $\mathcal{J}$ des structures presque complexes NK de la
sph{\`e}re $S^6$ est isomorphe {\`a} $\RM P(7)$.
\end{prop}
\begin{proof}
Le groupe d'isom{\'e}tries SO(7) agit transitivement sur les structures
presque complexes NK par
\begin{eqnarray*}
SO(S^6) \times \mathcal{J} & \to     & \mathcal{J} \\ 
(f,J)                      & \mapsto & f.J
\end{eqnarray*}
o{\`u} $\forall \, x \in M$
\[ (f.J)_x = (f_*)_x J_x (f_*)^{-1}_x
\]
Le groupe d'isotropie est isomorphe {\`a} $G_2$ car si $f \in SO(7)$
pr{\'e}serve $\omega$, il pr{\'e}serve aussi $\alpha$, d{\'e}finie par
(\ref{3forme}). Par cons{\'e}quent $\mathcal{J}$ est isomorphe {\`a}
$SO(7)/G_2 \simeq \RM P(7)$.
\end{proof}

\vs

Ceci ach{\`e}ve la classification des vari{\'e}t{\'e}s NK homog{\`e}nes
simplement connexes de dimension 6. En passant on a d{\'e}montr{\'e} la
conjecture \ref{conj} en dimension 6 et d'apr{\`e}s Nagy \cite{na2}, en toute
dimension. Les r{\'e}sultats sont r{\'e}sum{\'e}s dans les th{\'e}or{\`e}mes
\ref{theo} et \ref{theo2}.

\labelsep .5cm

\end{document}